\documentclass[onefignum,onetabnum]{siamonline190516}
\usepackage[utf8]{inputenc}
\usepackage{enumitem} 
\usepackage[normalem]{ulem}
\usepackage{amsmath}
\usepackage{amssymb}
\usepackage{subfig}
\usepackage{amsbsy,mathrsfs}
\usepackage{amsfonts}
\usepackage{hyperref}
\usepackage{tikz}
\usepackage{macros}
\usepackage{verbatim}
\usepackage{mathtools}
\usepackage{pgfplots}
\usepackage[margin=1in]{geometry}
\usepackage{xcolor,graphicx}
\usepackage{float}
\usepackage{bbm}
\usepackage{multirow}
\usepackage{multicol}
\usepackage{bigstrut}
\newcommand{\h}[1]{\mathbf{#1}}
\pgfplotsset{compat=1.17}

\headers{Implicit Regularization Effects}{Zhu, Hu, Lou and Yang}

\title{Implicit Regularization Effects of the Sobolev Norms in Image Processing\thanks{Submitted to the editors.
\funding{Y.~Yang was partially supported by NSF grant DMS-1913129. J.~Hu was partially supported by NSF CAREER grant  DMS-1654152. Y.~Lou was partially supported by NSF CAREER grant DMS-1846690. This paper is generated in the Summer Research Program for Women in Mathematics in Summer 2021. All authors acknowledge the generous support from the Mathematical Sciences Research Institute (MSRI). Y.~Yang acknowledges supports from Dr.~Max R\"ossler, the Walter Haefner Foundation and the ETH Z\"urich Foundation.  This work was done in part while Y.~Yang was visiting the Simons Institute for the Theory of Computing in Fall 2021.}}}

\author{
Bowen Zhu\thanks{New York University, New York, NY 10012, USA (\email{bz1010@nyu.edu})} 
\and 
Jingwei Hu\thanks{Department of Applied Mathematics, University of Washington, Seattle, WA 98195, USA (\email{hujw@uw.edu}).} \and  
Yifei Lou\thanks{Department of Mathematical Sciences, The University of Texas at Dallas, Richardson, TX 75080, USA (\email{yifei.lou@utdallas.edu})} 
\and 
Yunan Yang\thanks{Institute for Theoretical Studies, ETH Z\"urich, Z\"urich, 8092, Switzerland (\email{yunan.yang@eth-its.ethz.ch})}.
} 

\begin{document}

\maketitle

\begin{abstract}
In this paper, we propose to use the general $L^2$-based Sobolev norms, i.e., $H^s$ norms where $s\in \mathbb{R}$, to measure the data discrepancy due to noise in image processing tasks that are formulated as optimization problems. As opposed to a popular trend of developing regularization methods, we emphasize that an \textit{implicit} regularization effect can be achieved through the class of Sobolev norms as the data-fitting term. Specifically, we analyze that the implicit regularization comes from the weights that the $H^s$ norm imposes on different frequency contents of an underlying image. We further analyze the underlying noise assumption of using the Sobolev norm as the data-fitting term from a Bayesian perspective, build the connections with the Sobolev gradient-based methods and discuss the preconditioning effects on the convergence rate of the gradient descent algorithm, leading to a better understanding of functional spaces/metrics and the optimization process involved in image processing. Numerical results in full waveform inversion, image denoising and deblurring demonstrate the implicit regularization effects.

\end{abstract}

	
\section{Introduction}



 Digital images provide a powerful and intuitive way to represent the physical world. 
Unfortunately, noise is inevitable in the data that is taken or transmitted.  When recovering an underlying image from its corrupted measurements, one requires a \textit{fidelity} term to properly model the discrepancy of an imaging formation model as well as a \textit{regularization} term to refine the solution space of this inverse problem. The choice of such data fidelity term often depends on specific applications, specifically on the assumption of the noise distribution~\cite{bungert2020variational}. For example, a standard approach for additive Gaussian noise is the least-squares fitting. Using the maximum a posteriori (MAP) estimation, Aubert and Aujol \cite{aubert2008variational}
formulated a non-convex data fidelity term for multiplicative noise, which can be solved via a difference of convex algorithm \cite{li2016variational}.
In photon-counting devices such as x-ray computed tomography (CT) \cite{elbakri2002statistical,kak2001principles} and positron emission tomography (PET) \cite{vardi1985statistical}, the number of photons collected by a device follows a Poisson distribution, thus referred to as Poisson noise. Following the MAP of Poisson statistics,  the data discrepancy for Poisson noise can be modeled by a  log-likelihood form \cite{chowdhury2020non,chowdhury2020poisson,le2007variational}. Since the nonlinearity of such data fidelity causes computational difficulties, a popular approach in CT reconstruction adopts a weighted least-squares model \cite{thibault2007three} as the data fitting term.

 To date, major research interests in image processing community have  focused on developing regularization methods by exploiting the prior knowledge and/or the special structures of an imaging problem. For instance, the classic Tikhonov regularization~\cite{tikhonov1943stability} returns a smooth output in an attempt to remove the noise, however, at the cost of smearing out important structures and edges. Total variation (TV)  \cite{rudin1992nonlinear} is an edge-preserving regularization in that it tends to diffuse along the edges, rather than across, but TV causes a staircasing (blocky) artifact. As remedies, 
  total generalized variation (TGV)  \cite{bredies2010total} and fractional-order TV (FOTV) \cite{zhang2015total} were proposed to preserve higher-order smoothness.  In addition, non-local regularizations  \cite{lou2010image,zhang2010bregmanized} based on patch similarities \cite{buades2005review} work well for textures and repetitive patterns in an image.

Instead of proposing explicit regularization models, we reveal in this paper that implicit regularization effects can be achieved by using only the $L^2$-based Sobolev norms as a data fidelity term. Recall that a  Sobolev space is a vector space of functions equipped with a norm that combines the $L^p$ norms of the function and its derivatives up to a given order. We are particularly interested in the $L^2$-based Sobolev spaces, often referred to as the $H^s$ spaces for $s\in \mathbb R,$  since they are well-studied and widely used. Note that an $H^s$ space is also a Hilbert space with a well-defined inner product. Its associated norm, which we refer to as the $H^s$ norm,  is naturally equipped with a particular form of weighting in the Fourier domain. Both the order of biasing (e.g., towards either low or high frequencies) and the strength of biasing can be controlled by the choice of $s\in\mathbb R$. When $s=0$, it reduces to the standard $L^2$ norm with equal weights on all the frequencies due to Parseval's identity. Since $H^s$ is a generalization of the $L^2$ norm, using the $H^s$ norm undoubtedly leads to improved results when the parameter $s$ is appropriately chosen according to the prior information, e.g., noise spectra.
On the other hand, the $H^s$ norms offer additional flexibility by choosing  $s$ to achieve either smoothing ($s<0$) or sharpening ($s>0$) effects depending on the noise type in an input image. 
It was analyzed in~\cite{engquist2020quadratic} that the class of the $H^s$ norms brings a preconditioning effect as an objective function, thus altering the stability of the original inverse problem. In~\cite{yang2020anderson}, a particular frequency bias of the $H^s$ norm was utilized to accelerate fixed-point iterations when seeking numerical solutions to elliptic partial differential equation
(PDEs).


The introduction of Sobolev spaces was significant for the development of functional analysis~\cite{sobolev1963applications}
and various applications related to PDEs~\cite{evans98} such as the finite element method~\cite{szabo1991finite}.
There have been relevant works to the Sobolev norms in image processing and inverse problems. For example, the $H^{-1}$ semi-norm is closely related to the quadratic Wasserstein ($W_2$) metric from optimal transportation~\cite{villani2003topics} under both the asymptotic regime~\cite{otto2000generalization} and the non-asymptotic regime~\cite{peyre2018comparison}. This connection has been utilized in many applications~\cite{engquist2020quadratic,papadakis2014optimal} such as Bayesian inverse problems~\cite{dunlop2020stability}. 
Another close connection comes from works on the Sobolev gradient~\cite{neuberger2009sobolev}, in which the gradient of a given functional is taken with respect to the inner product induced by the underlying Sobolev norm~\cite{calderMY10,sundaramoorthi2007sobolev} with demonstrated effects in sharpening and edge-preserving. 

In this paper, we illustrate the implicit regularization effects of the $H^s$ norm as a data fitting term on a toy example of deblurring a square image, together with two geophysical applications of image denoising and full waveform inversion. In those examples, we use only the $H^s$ norm as a data fidelity term in the objective function without any regularization term. The final reconstructions mitigate the impact of the noise, reflecting the implicit regularization effects. This approach is particularly effective when the spectral contents of the noise are well-separated from the spectral contents of the actual image.
Since some natural images have a broad bandwidth with spectral contents spreading out in the frequency domain, the implicit regularization by $H^s$ alone may not effectively preserve the important features. In those scenarios, it is beneficial to incorporate, for example, the total variation as a regularization term together with the $H^s$ norm as the data fidelity. We acknowledge that using the $H^s$ norm as the data fidelity term together with the total variation regularization has been intensively studied in~\cite{osher2003image,lieu2008image}. In this work, we generalize their approaches by considering $s$ as a tunable hyperparameter in practical implementations and proposing a more efficient algorithm by the alternating direction method of multipliers (ADMM) \cite{boyd2011distributed,glowinski1975approximation}.

The main contributions of this work are threefold. First, we propose to use the $H^s$ norms as a novel data-fitting term to effectively utilize their implicit regularization effects for noise removal. Second, we analyze the underlying noise assumption of using the $H^s$ norms as the objective function from a Bayesian perspective, its connections to the Sobolev gradient flow, and the resulting preconditioning effects on the convergence rate. Such analysis contributes to a better understanding of the advantages by using the $L^2$-based Sobolev norms in image processing. Lastly, we present a series of computational approaches to calculating the $H^s$ norms under different setups.



The rest of the paper is organized as follows. \Cref{sect:analysis} devotes to the analysis of the Sobolev norms, including the implicit regularization effects, the noise assumption from a Bayesian perspective, the connections to the $W_2$ distance, the Sobolev gradient, and the preconditioning effects. We describe three approaches for computing the $\cH^s$ norm in \Cref{sec:Hs_dist} under different boundary conditions and choices of $s$. In \Cref{sect:exp}, we conduct experiments on geographical examples to demonstrate different scenarios where the weak norm ($s<0$) and the strong norm ($s>0$) are preferred, respectively. \Cref{sect:TVHs} revisits the $H^s$+TV model~\cite{osher2003image,lieu2008image} with a tunable parameter $s$ and an efficient algorithm for image deblurring. 
Conclusions follow in \Cref{sect:conclusion}. %

\section{Analysis on Sobolev Norms}
\label{sect:analysis}

In this section, we   briefly review the definitions and properties of the  $L^2$-based Sobolev norms, followed by discussing the implicit  regularization effects in \Cref{sec:Hs_analysis}.
We draw connections of the Sobolev norms to a Bayesian interpretation of data fidelity in~\Cref{sec:baysian}, the quadratic Wasserstein distance~\cite{villani2003topics} in~\Cref{sect:W2}, and the Sobolev gradient~\cite{calderMY10} in~\Cref{sect:sob_grad}. Lastly in~\Cref{sect:conv-rate}, we discuss how the choice of the Sobolev norm can affect the convergence rate of the gradient descent algorithm. 

\subsection{$\H^s$ Sobolev Space} \label{sec:Hs_def}
There are two common ways to define the $L^2$-based Sobolev norm. One is based on the Sobolev space $W^{k,p}(\R^d)$ for a nonnegative integer $k$; see \cref{def:Wkp_def}.   
\begin{defn}[Sobolev Space $W^{k,p}(\R^d)$\label{def:Wkp_def}] Let $1\leq p<\infty$ and $k$ be a nonnegative integer. If a function $f$ and its weak derivatives $D^{\alpha}f=\frac{\partial^{|\alpha|}f}{\partial x_1^{\alpha_1}\cdots \partial x_d^{\alpha_d}}$, $|\alpha|\leq k$ all lie in $L^p(\R^d)$, where $\alpha$ is a multi-index and $|\alpha|=\sum_{i=1}^d \alpha_i$, we say $f\in W^{k,p}(\R^d)$ and define the $W^{k,p}(\R^d)$ norm of $f$ as
\begin{equation}\label{eq:Wkp_def}
    \|f\|_{W^{k,p}(\R^d)} := \left(\sum_{|\alpha|\leq k} \| D^{\alpha} f \|_{L^p(\R^d)}^p\right)^{1/p}.
\end{equation}
\end{defn}
In this work we focus on the $L^2$-based Sobolev space $W^{k,2}$, which is a Hilbert space.

 While~\Cref{def:Wkp_def} is concerned with integer-regularity spaces, there exists a natural extension to a more general $L^2$-based Sobolev space $W^{s,2}(\R^d)$ for an arbitrary scalar $s\in\R$ through the Fourier transform. This leads to the second definition of the Sobolev space. Specifically, we define 
\begin{eqnarray}
\F f(\xi) = \hat{f}(\xi) = (2\pi)^{-\frac{d}{2}} \int_{\mathbb{R}^d} f(x) e^{-ix\cdot \xi}dx,
\end{eqnarray}
where $\F$  denotes the Fourier transform. We further denote $\F^{-1}$ as the inverse Fourier transform,  $I$ as the identity operator, $\aver{\xi} := \sqrt{1 + |\xi|^2}$, and $\mathcal{S}^\prime(\R^d)$ as the space of tempered distributions.

\begin{defn}[Sobolev Space $\cH^{s}(\R^d)$\label{def:Hs_def}] Let $s \in \R$, the Sobolev space $\cH^{s}$ over $\R^d$ is given by
\begin{align}
    \cH^{s}(\R^d) := \left\{ f \in \mathcal{S}'(\R^d) : \F^{-1} \left[ \aver{\xi}^{s} \F f \right] \in L^2(\R^d) \right\}.
\end{align}
The space $ \cH^{s}(\R^d)$ is equipped with the norm
\begin{align} \label{eq:Hs_def}
    \|f\|_{\cH^{s}(\R^d)} := \left\| \F^{-1} \left[ \aver{\xi}^{s} \F f \right] \right\|_{L^2(\R^d)} = \left\| \mathcal{P}_sf \right\|_{L^2(\R^d)},
\end{align}
where the operator $\mathcal{P}_s:=(I-\Delta)^{s/2}$.
\end{defn}

When $s=0$, the $H^{s}(\R^d)$ space (norm) reduces to the standard $L^2$ space (norm). 
One can show that $W^{k,2}(\R^d) = H^k(\R^d)$ for any 
integer $k$~\cite{arbogastmethods}. We remark that  $ \|f\|_{H^{k}(\R^d)} \neq    \|f\|_{W^{k,2}(\R^d)} $ for the same $k$ in general, but the two norms are equivalent, which can be shown through Fourier transforms. Hereafter, we  mainly focus on  $H^{s}(\R^d)$ for $s\in \R$, due to its better generality. 

\subsection{Implicit Regularization Effects of the $H^s$ Norms}
\label{sec:Hs_analysis}
Without loss of generality,
we consider the following data formation model based on a linear inverse problem,
\begin{equation} \label{EQ:Lin IP}
	 f_\sigma = \mathcal{A} u + n_\sigma,
\end{equation}
where $f_\sigma$ denotes the noisy measurements with an additive Gaussian noise $n_\sigma$ of standard deviation $\sigma$, and $\mathcal A$ denotes a linear degradation operator. A general inverse problem is posted as recovering an underlying image $u$ from the data $f_\sigma$ with the knowledge of $\mathcal A.$ If $\mathcal A$ is the identity operator, i.e., $\mathcal A = I,$ this problem is referred to as \textit{denoising}. If $\mathcal A$ can be formulated as a convolution operator with a blurring kernel, it is called image \textit{deblurring} or \textit{deconvolution.} 


We assume the linear operator $\mathcal{A}$ is asymptotically diagonal in the Fourier domain such that
\begin{equation}\label{EQ:IP Symbol}
	\widehat{\mathcal{A} u}(\xi)\sim \aver{\xi}^{-\alpha} \hat{u}(\xi),
\end{equation}
where $\alpha\in\R$, the hat symbol  denotes the  Fourier transform with frequency coordinate $\xi$, and $\sim$ refers to the relationship that both sides are asymptotically on the same order of magnitude. 
When $\alpha> 0$, we say the operator $\mathcal{A}$ is ``smoothing''. The value of $\alpha$ can describe  to some extent the degree of ill-conditionedness (or difficulty) of solving an inverse problem~\cite{bal2012introduction} in the sense that the larger the $\alpha$ is, the more ill-posed the associated inverse problem becomes. 

We examine the regularization effects of using the $\cH^s$ norm  defined in~\eqref{eq:Hs_def} to quantify the data misfit. In other words, we seek a solution of the inverse problem \eqref{EQ:Lin IP} by minimizing
\begin{equation}\label{eq:Hs_obj_old}
	\Phi_{\cH^s}(u) := \frac{1}{2}\|\mathcal{A}u-f_\sigma\|^2_{\cH^s}=\frac{1}{2}\|\mathcal{P}_s(\mathcal{A}u-f_\sigma)\|^2_{L^2}=\frac{1}{2}\int_{\R^d}\aver{\xi}^{2s}|\widehat{\mathcal{A}u}(\xi) - \widehat{f_\sigma}(\xi)|^2 d\xi,
\end{equation}
without any additional regularization term. The minimizer of $\Phi_{\cH^s}(u)$ has a closed-form solution, i.e.,
\begin{equation}\label{EQ:Hs Inversion Phy}
	u =\Big(\mathcal{A}^* \mathcal{P}_s^*\mathcal{P}_s \mathcal{A} \Big)^{-1} \mathcal{A}^* \mathcal{P}_s^*\mathcal{P}_s f_\sigma,
\end{equation}
where $\mathcal{A}^*$ is the adjoint operator of $\mathcal{A}$ under the $L^2$ inner product and  $\mathcal{P}_s= (I-\Delta)^{s/2}.$ 
Note that $\mathcal{P}_s^* = \mathcal{P}_s$ as $\mathcal{P}_s$ is self-adjoint. By comparing \eqref{EQ:Hs Inversion Phy} with the standard least-squares solution, we conclude that the $\cH^s$-based inversion  can be seen as a weighted least-squares method if $s\neq 0$. 
\begin{rmk}
A variant of~\eqref{eq:Hs_obj_old} is to use the $\dot{H}^s$ semi-norm instead of the standard $H^s$ norm. That is, we replace $\aver{\xi}^{2s} = (1 + |\xi|^2)^s$ by $|\xi|^{2s}$, and the objective function becomes
\begin{equation}\label{eq:Hs_obj_semi}
	 \Phi_{\dot{\cH}^s}(u)= \frac{1}{2}\|\mathcal{A}u-f_\sigma\|^2_{\dot{\cH}^s}:=\frac{1}{2}\int_{\R^d} |\xi|^{2s}|\widehat{\mathcal{A}u}(\xi) - \widehat{f_\sigma}(\xi)|^2 d\xi.
\end{equation}
The  frequency bias from $\Phi_{\dot{\cH}^s}$ is more straightforward to analyze than the one from  $\Phi_{\cH^s}(u)$, as the weight in front of each frequency is precisely an algebraic factor $|\xi|^{s}$. If $f\in \cH^{s}$ for $s>0$, we  have $||f||_{\dot{\cH}^s} < \infty$. However, this is not the case for $s<0$. For example, a function $f$ may have a finite $\cH^{-1}$ norm, but if $\int f dx \neq 0$, it does not have a well-defined $\dot{\cH}^{-1}$ norm.
\end{rmk}
\begin{rmk}
If $s_1,s_2\in \R$ and $s_1 < s_2$, then $\cH^{s_2} \subset \cH^{s_1}$ is continuously embedded. In other words, we specify the order among all $H^s$ spaces, e.g., $H^2 \subset H^1 \subset L^2\subset H^{-1} \subset H^{-2}$.
\end{rmk}

We consider  the following three  scenarios to illustrate the implicit regularization effects of $\Phi_{\cH^s}$ as an objective function. A similar analysis applies to $\Phi_{\dot{\cH}^s}$. 
\begin{itemize}
    \item When $s=0$, the solution \eqref{EQ:Hs Inversion Phy} reduces to the standard least-squares solution, i.e., $u= \mathcal{A}^\dagger f_\sigma,$ where $\mathcal{A}^\dagger$ is the Moore--Penrose inverse
operator of $\mathcal{A}$. Without any regularization term, this solution inevitably overfits the noise in the observation $f_\sigma$. 

\item When $s>0$, $\mathcal{P}_s$ can be regarded as a differential operator, which amplifies high-frequency contents of $f_\sigma$. If the noise in $f_\sigma$ is also high-frequency, the overfitting phenomenon caused by $\mathcal{P}_s$ is even worse than the standard least-squares solution. On the other hand, if $f_\sigma$ is corrupted by lower-frequency noise, the weighted least-squares would avoid overfitting. 

\item When $s<0$,  $\mathcal{P}_s$ is an integral operator, meaning that applying  $\mathcal{P}_s$ to $f_\sigma$ suppresses high-frequency components. The noisy content in $f_\sigma$ does not fully ``propagate'' into the reconstructed solution $u$. The inverse problem is less sensitive to the high-frequency noise in $f_\sigma$, indicating the improved well-posedness. Again, this property becomes disadvantageous if $f_\sigma$ is subject to lower-frequency noise. 
\end{itemize}
Based on the above three different types of scenarios, it is clear that the $H^s$ norm causes a particular weight on the frequency contents of the input function according to the choice of $s$. We will later refer to this property as the \textit{spectral bias} of the $H^s$ norm.

\begin{rmk}\label{rmk:goodbad}
To summarize, if the data is polluted with high-frequency noise, using a weak norm as the objective function alone  improves the posedness of the   inverse data-fitting problem without the help of any regularization term. On the other hand,  a potential disadvantage of the weaker norm is that the objective function not only implicitly suppresses the higher-frequency noisy content but also the higher-frequency component of the noise-free data. Consequently, the reconstruction loses the high-frequency resolution, as illustrated in~\cite[Figure 4]{engquist2020quadratic}.
\end{rmk}

\begin{rmk}
One can also generalize~\eqref{EQ:Lin IP} to a nonlinear inverse problem. The main properties of the $H^s$ norm will remain, but the analysis would be less straightforward. In~\Cref{sect:exp_fwi}, we present such a nonlinear example and numerically demonstrate the benefits of using the $H^s$ norm.
\end{rmk}

Next, we demonstrate the aforementioned properties regarding $s=0$, $s>0$ and $s<0$ through  numerical examples of reconstructing a (discrete) image $u$ from \eqref{EQ:Lin IP}
by minimizing the discretized objective function $$\Phi_{\cH^s}(u) = \frac{1}{2}\|P_s(Au-f_\sigma)\|^2_{L^2},$$ where $P_s$ is a proper discretization of the continuous operator $\mathcal{P}_s$, and $A$ denotes the linear operator $\mathcal{A}$ in the matrix form; please refer to Section~\ref{sec:Hs_dist} for  discretization details. Applying the gradient descent algorithm with a fixed step size $\eta$ to minimize the objective function $\Phi_{\cH^s}(u)$ yields the following iterative step:
\begin{equation}\label{eq:Hs_GD}
    u^{(n+1)}  = u^{(n)} - \eta \nabla \Phi(u^{(n)}) =  u^{(n)} - \eta  A^TP_s^T P_s (Au^{(n)}-f_\sigma).
\end{equation}


\begin{figure}
\centering
\subfloat[Blurry Input]{\includegraphics[width = 0.15\textwidth]{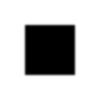}\label{fig:square_no_noise_input}}
\hspace{0.1cm}
\subfloat[$s=1$]{\includegraphics[width = 0.15\textwidth]{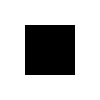}}
\hspace{0.1cm}
\subfloat[$s=0.5$]{\includegraphics[width = 0.15\textwidth]{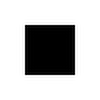}}
\hspace{0.1cm}
\subfloat[$s=0$]{\includegraphics[width = 0.15\textwidth]{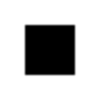}}
\hspace{0.1cm}
\subfloat[$s=-0.5$]{\includegraphics[width = 0.15\textwidth]{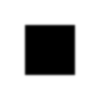}}
\hspace{0.1cm}
\subfloat[$s=-1$]{\includegraphics[width = 0.15\textwidth]{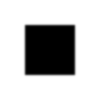}}
\caption{Effects of minimizing $\Phi_{\cH^s}$  with different choices of $s$. The reconstructed solutions gradually transition from sharp to blurry after the same number of gradient descent iterations, showing that strong norms ($s>0$) are better at sharpening. \label{fig:square_deblur_nonoise_different_s}}
\end{figure}


We apply~\eqref{eq:Hs_GD} to a simple example of image deblurring. Consider a binary image of size $100\times 100$ with a black square in the middle to be the ground-truth, referred to as the Square image. The linear operator $A$ can be formulated as a convolution with $15\times 15$ Gaussian  kernel of standard deviation $1$, which can be implemented through {\tt fspecial(`gaussian',15,1)} in Matlab. The blurry image is further corrupted by an additive Gaussian noise with standard deviation $\sigma$. 

When $\sigma=0$, the input image is  blurry but not noisy, as seen in~\Cref{fig:square_no_noise_input}. We  show reconstructed images by minimizing  $\Phi_{\cH^s}$  with different choices of $s$ via~\eqref{eq:Hs_GD}. The five values of $s$ in~\Cref{fig:square_deblur_nonoise_different_s} cover all scenarios: $s=0$, $s>0$ and $s<0$. After running 
$100$ iterations of the gradient descent algorithm~\eqref{eq:Hs_GD} with the same step size $\eta=1$, we observe  in \Cref{fig:square_deblur_nonoise_different_s} a gradual transition  from sharp to blurry reconstruction results as $s$  decreases from $s=1$ to $s=-1$. This is  aligned with our earlier discussion that
the operator $\mathcal P_s$ for positive $s$ is a differential operator, which boosts the higher-frequency content of $A^T (Au^{(n)}-f_\sigma)$, which is the gradient when the $L^2$ norm becomes the objective function. Consequently, it accelerates the gradient descent algorithm to converge to the sharp ground truth, as the only missing information in the blurry input is precisely in the high-frequency domain. In summary, strong norms ($s>0$) are good at sharpening.

\begin{figure}
\centering
\subfloat[Noisy Input]{\includegraphics[width = 0.15\textwidth]{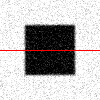}\label{fig:square-blur01}}
\hspace{0.1cm}
\subfloat[$s=-1$]{\includegraphics[width = 0.15\textwidth]{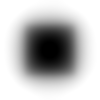}}
\hspace{0.1cm}
\subfloat[$s=-0.75$]{\includegraphics[width = 0.15\textwidth]{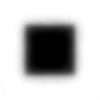}}
\hspace{0.1cm}
\subfloat[$s=-0.5$]{\includegraphics[width = 0.15\textwidth]{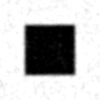}}
\hspace{0.1cm}
\subfloat[$s=-0.25$]{\includegraphics[width = 0.15\textwidth]{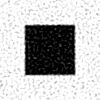}}
\hspace{0.1cm}
\subfloat[$s=0$ ($L^2$)]{\includegraphics[width = 0.15\textwidth]{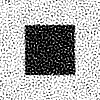}}\\
\subfloat[Noisy Input]{\includegraphics[width = 0.15\textwidth]{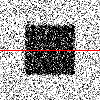}\label{fig:square-blur05}}
\hspace{0.1cm}
\subfloat[$s=-1$]{\includegraphics[width = 0.15\textwidth]{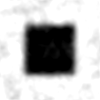}}
\hspace{0.1cm}
\subfloat[$s=-0.75$]{\includegraphics[width = 0.15\textwidth]{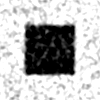}}
\hspace{0.1cm}
\subfloat[$s=-0.5$]{\includegraphics[width = 0.15\textwidth]{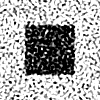}}
\hspace{0.1cm}
\subfloat[$s=-0.25$]{\includegraphics[width = 0.15\textwidth]{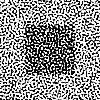}}
\hspace{0.1cm}
\subfloat[$s=0$ ($L^2$)]{\includegraphics[width = 0.15\textwidth]{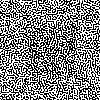}}
\caption{Deblurring the Square image by minimizing $\Phi_{\cH^s}(u)$. The top row presents the blurry noisy input with $\sigma = 0.1$ and reconstruction results of different $s$ values. A noisier case ($\sigma = 0.5$) is illustrated in the bottom row. \label{fig:square_deblur}}
\end{figure}

We then examine the influence of noise on the reconstructions by minimizing  the $\Phi_{\cH^s}$ functional. For this purpose, we add different amounts of noises, i.e., $\sigma = 0.1$ and $\sigma=0.5$, to the same blurry image (shown in~\Cref{fig:square_no_noise_input}), leading to noisy and blurry data shown in~\Cref{fig:square-blur01} and~\Cref{fig:square-blur05}, respectively.  Again, we reconstruct the images by running 100 iterations of gradient descent with the same step size. The top row of~\Cref{fig:square_deblur} corresponds to a smaller noise level ($\sigma = 0.1$). The $L^2$-based method, i.e., $s=0,$ clearly suffers from  overfitting the noise, as the reconstruction is even noisier than the input. The best result is achieved at $s= -0.5,$ while the reconstructed images are over smooth as $s$ decreases. This set of tests shows both  advantages and  potential limitations of weak norms  ($s<0$) as addressed in~\Cref{rmk:goodbad}. The bottom row of~\Cref{fig:square_deblur} corresponds to a larger noise level ($\sigma = 0.5$), when the overfitting phenomenon is more severe not only for the $L^2$ norm, but also for the cases of $s= -0.5$ and $s=-0.25$. The best reconstruction occurs at $s=-1$, where the spectral bias of the objective function towards lower-frequency contents of the residual (the difference between the current iterate and the input image)
is the strongest. That is, the weighting coefficients on the low-frequency components are much bigger in contrast to the ones on the high-frequency ones due to the rapid decay of function $\langle \xi\rangle^{-1}$ compared to $\langle \xi \rangle ^{-0.5}$. The comparison between two noise levels also implies that the best choice of $s$ is data-dependent. One heuristic principle is that the noisier the input is, the weaker objective function (smaller $s$) one should choose to avoid overfitting the noise.

\begin{figure}
    \centering
    \subfloat[$\sigma=0.1$]{\includegraphics[width = 1.0\textwidth]{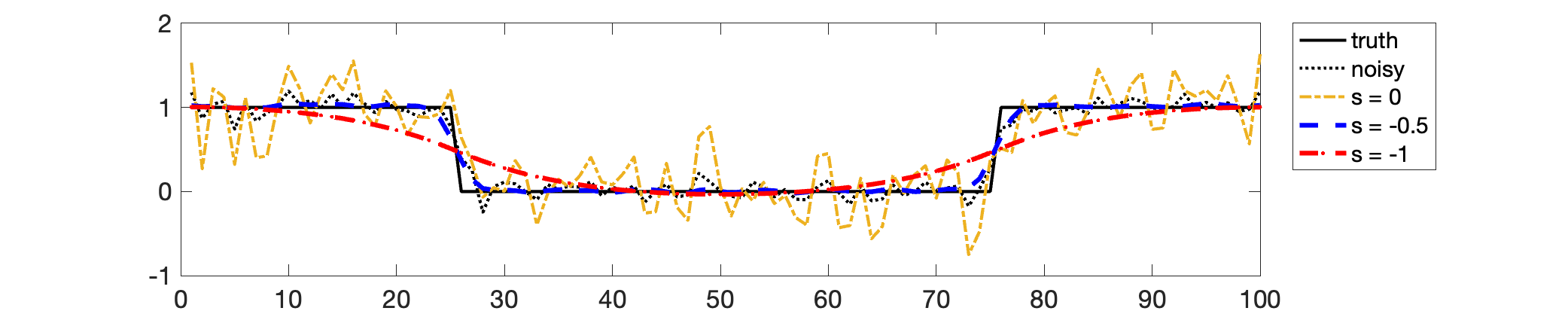}\label{fig:zoomin_01}}\\
    \subfloat[$\sigma=0.5$]{\includegraphics[width = 1.0\textwidth]{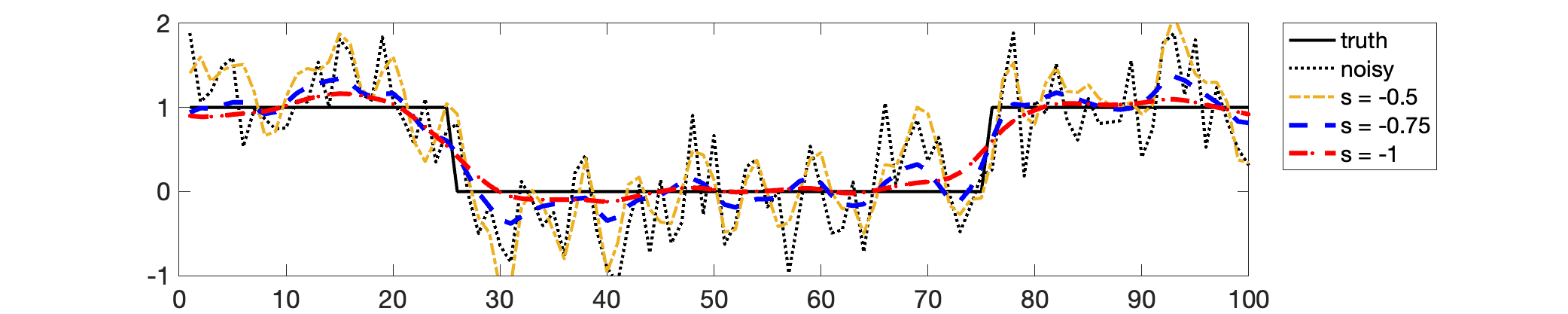}\label{fig:zoomin_05}}    
    \caption{The zoom-in view for different choice of $s$ at the cross section (the red line) illustrated in~\Cref{fig:square-blur01} and~\Cref{fig:square-blur05}, respectively.}
    \label{fig:square_deblur_zoom_in}
\end{figure}

In~\Cref{fig:square_deblur_zoom_in}, we show the cross-sections of   2D images; the location of the cross-section is indicated by the red lines in~\Cref{fig:square-blur01} and~\Cref{fig:square-blur05}. In~\Cref{fig:zoomin_01}, the 1D plots clearly show the over-smoothing artifact for $s=-1,$ and the construction of $s=-0.5$ is closest to the ground truth. In contrast, the case $s=-0.5$ is no longer good enough to ``smooth'' out the stronger noise
in~\Cref{fig:zoomin_05},  and the result from $s=-1$ turns out to be the best fit.

\subsection{A Bayesian Interpretation}\label{sec:baysian}
The choice of the data fidelity term in image processing can be derived from a Bayesian approach under a proper assumption on the noise distribution of the data~\cite{bungert2020variational}. In this subsection, we present the noise assumption associated with the proposed $H^s$ data fidelity term~\eqref{eq:Hs_obj_old} under the Bayesian framework.

One major advantage of the Bayesian approach is to account for the uncertainty in the data which will be propagated to the solution to the inverse problem. It combines a probabilistic model for the observed data $f_\sigma$ with a density function $\BP(f_\sigma|u)$ and a probability distribution $\BP(u)$ representing the prior knowledge regarding the unknown $u$. Bayes' theorem provides a way to construct the posterior distribution, denoted as $\BP( u| f_\sigma)$, where 
\begin{align}
\label{eq:bayes_formal}
\BP(u|f_\sigma) = \frac{\BP(f_\sigma|u)\BP( u)}{\BP(f_\sigma)}.
\end{align}
The posterior distribution $\BP( u|f_\sigma)$ can be regarded as the solution to the Bayesian inverse problem, which differs from the deterministic framework of solving inverse problems that returns a single value of $u$, e.g., the minimizer of~\eqref{eq:Hs_obj_old}. 

Although the Bayesian and the deterministic approaches are quite different, there are connections when we try to find the maximum a posteriori (MAP) estimation. Without loss of generality, we consider that the prior distribution $\BP( u)$ follows the normal distribution $\mathcal{N}(0,C)$ and $C$ is invertible. Then maximizing the posterior distribution $\BP( u|f_\sigma)$ is equivalent to the following minimization problem~\cite[Sec.~4.3]{Dashti2017}, 
\begin{equation} ~\label{eq:bayes_misfit}
u^* = \argmin_u \mathcal{E}(u;f_\sigma) + \frac{1}{2}\langle u,C^{-1} u \rangle_{L^2},
\end{equation}
where $\mathcal{E}(u; f_\sigma) = -\log \BP(f_\sigma|u)$ is commonly known as the negative log-likelihood function. 
Consider the inverse problem model~\eqref{EQ:Lin IP} where we assume the additive noise $n_\sigma\sim \mathcal{N}(0, \Gamma)$. We then have
\begin{equation*} 
\BP(f_\sigma|u) = \mathcal{N}(\mathcal{A}u, \Gamma) \propto \exp \left(  - \frac{1}{2}  \|\mathcal{A}u - f_\sigma \|^2_\Gamma \right), \quad \langle\cdot,\cdot\rangle_\Gamma := \langle\cdot,\Gamma^{-1}\cdot\rangle_{L^2},
\end{equation*}
after ignoring the normalizing constant. If the inverse covariance operator $\Gamma^{-1} = \mathcal{P}_s^* \mathcal{P}_s$ for $\mathcal{P}_s$ defined in~\Cref{sec:Hs_analysis}, we then have
\[
\mathcal{E} (u;f_\sigma) = -\log \BP(f_\sigma|u)  \propto \frac{1}{2} \|\mathcal{A}u - f_\sigma \|^2_\Gamma = \frac{1}{2} \|\mathcal{P}_s \left(  \mathcal{A}u - f_\sigma \right)\|^2_{L^2}, 
\]
which reduces to our $H^s$ objective function~\eqref{eq:Hs_obj_old}. 

Based on the above Bayesian interpretation, using~\eqref{eq:Hs_obj_old} as the data-fidelity term is equivalent to a data noise assumption $n_\sigma\sim \mathcal{N}(0, \left(\mathcal{P}_s^* \mathcal{P}_s\right)^{-1})$ in the Bayesian  framework. Note that the $L^2$ norm corresponds to $n_\sigma\sim \mathcal{N}(0, I)$, the standard Gaussian. This perspective again demonstrates that we can enforce prior information to achieve implicit regularization effects through the data fidelity (likelihood function) term. For example, $n_\sigma\sim \mathcal{N}(0, (I-\Delta)^{-1})$ ($s = 1$) supposes a smooth additive noise while $n_\sigma\sim \mathcal{N}(0, I-\Delta)$ ($s=-1$) assumes that the noise lacks of smoothness. This interpretation also extends to the seminorm $\Phi_{\dot{H}^s}$~\eqref{eq:Hs_obj_semi}, which corresponds to $n_\sigma\sim \mathcal{N}(0, (-\Delta)^{-s})$.

\subsection{Relationship with the $W_2$ Distance}\label{sect:W2}
Here, we review  a connection between the Sobolev norms and the quadratic Wasserstein ($W_2$) distance~\cite{villani2003topics}  to provide a better understanding of both metrics.
The Wasserstein distance defined below is associated to the cost function $c(x,y) = |x-y|^p$ in the optimal transportation problem. 
\begin{definition}[Wasserstein Distance]
We denote by $\mathscr{P}_p(\Omega)$ the set of probability measures with finite moments of order $p$. For $1\leq p<\infty$,  
\begin{equation}\label{eq:static}
W_p(\mu,\nu)=\left( \inf _{T_{\mu,\nu}\in \mathcal{M}}\int_{\Omega}\left|x-T_{\mu,\nu}(x)\right|^p d\mu(x)\right) ^{\frac{1}{p}},\quad \mu, \nu \in \mathscr{P}_p(\Omega),
\end{equation}
where $\mathcal{M}$ is the set of all maps that push forward $\mu$ into $\nu$. Note that $W_2$ corresponds to the case $p=2$.
\end{definition}
An asymptotic connection between the $W_2$ metric and the $\cH^s$ norm was first provided in~\cite{otto2000generalization} given the two probability distributions under comparison are close enough such that the linearization error is small. Consider $\mu$ as the probability measure and $d\pi$ as an infinitesimal perturbation that has zero total mass. Then
\begin{equation}\label{EQ:W2-Hm1 Asym}
	W_2(\mu, \mu+d\pi)=\|d\pi\|_{\dot{\cH}_{(d\mu)}^{-1}}+\smallO (d\pi).
\end{equation}
We remark that $\dot{\cH}_{(d\mu)}^{-1}$ is the weighted $\dot{\cH}^{-1}$ semi-norm. We refer readers to~\cite[Sec.~7.6]{villani2003topics} for its detailed definition.

A connection between $W_2$ and $\dot{\cH}^{-1}$ under a non-asymptotic regime was later presented in~\cite{peyre2018comparison}. If both $f dx=d\mu$ and $g dx=d\nu$ are bounded from below and above by constants $c_1$ and $c_2$,  we have the following \emph{non-asymptotic} equivalence between $W_2$ and $\dot{\cH}^{-1}$~\cite{peyre2018comparison}, 
\begin{equation}\label{EQ:W2-Hm1}
	\frac{1}{c_2} \|f-g\|_{\dot {\cH}^{-1}} \le W_2(\mu, \nu) \le \frac{1}{c_1} \|f-g\|_{\dot {\cH}^{-1}}.
\end{equation}

Note that in both the asymptotic and the non-asymptotic regimes, the $W_2$ metric shares a similar spectral bias as the $\dot {\cH}^{-1}$ semi-norm, up to a weighting function. Thus, the implicit regularization properties for the case $s=-1$ discussed in~\Cref{sec:Hs_analysis} can extend to the quadratic Wasserstein metric. This finding explains the improved stability of the Wasserstein metric in inverse problems from various applied fields, including machine learning~\cite{arjovsky2017wasserstein}, parameter identification~\cite{yang2021optimal}, and full-waveform inversion~\cite{yang2018application}.

\subsection{Relationship with the Sobolev Gradient Flow} \label{sect:sob_grad}
The well-known heat equation $u_t = \Delta u$ where  $u:\Omega  \mapsto
\mathbb{R}$ ($\Omega$ is an open subset of
$ \mathbb{R}^2$ with smooth boundary $\partial \Omega$) can be seen as the gradient flow of the energy functional
\[
E(u) = \frac 1 2 \int_\Omega |\nabla u|^2 dx = \frac 1 2\|\nabla
u\|_{L^2}^2 ,
\]
with respect to the $L^2$ inner product $\langle v, w
  \rangle_{L^2}  = \int_{\Omega} v\, w\,dx$.  A different gradient
flow can be derived from a more general inner product, for example, based on the Hilbert space $\cH^{s}$ in~\Cref{def:Hs_def} for any $s\in\R$. An inner
product on the Sobolev space $\cH^1(\Omega)$~\cite{evans98,sundaramoorthi2007sobolev} can be
defined as 
$$
g_\lambda(v, w) = (1-\lambda) \langle v, w
  \rangle_{L^2} + \lambda \langle v, w \rangle_{\cH^1} = \langle v, w
  \rangle_{L^2} + \lambda \langle v, w \rangle_{\dot{\cH}^1} ,
$$ 
for any $\lambda >0$ and $\langle v, w \rangle_{\dot{\cH}^1}  = \langle \nabla v, \nabla w \rangle_{L^2}$. If we are only interested in periodic functions   on the domain $\Omega$, the gradient 
operators considered here are equipped with the periodic boundary condition. When $\lambda=0$, $g_\lambda(v, w)$ reduces the conventional $L^2$ inner product, and when $\lambda = 1$, it becomes the standard $\cH^1$ inner product: $\langle v,w \rangle_{{\cH}^1}  = \langle  v,  w \rangle_{L^2} + \langle \nabla v, \nabla w \rangle_{L^2}$. Calder \emph{et al.}~\cite{calderMY10} exploited  a general Sobolev gradient flow for image processing and established the well-posedness of the  Sobolev gradient flow $ u_t = (I-\lambda\Delta)^{-1}\Delta u$  in both the forward and the backward directions of minimizing $E(u)$. Specifically worth noticing is that the backward direction can be regarded as a sharpening operator \cite{liu2020image,lou2013video}. 

Without loss of generality, we set $\lambda=1$ when studying a connection  between the Sobolev gradient and the gradient of the $\cH^s$ norm as the energy functional.
Given any energy (objective) functional $E(u)$, an inner product based on the Sobolev metric $\cH^1(\Omega)$ gives a specific gradient formula
\begin{equation}
      \nabla_{\cH^1} E(u)  =(I-\Delta)^{-1}       \nabla_{L^2} E(u),
\end{equation}
such that
\begin{equation}
   \langle     \nabla_{\cH^1} E(u) , v \rangle_{\cH^1}  =   \langle \nabla_{L^2} E(u) , v \rangle_{L^2} = \lim_{\epsilon\rightarrow 0} \frac{E(u+\epsilon v) - E(u)}{\epsilon},\quad \forall v\in \cH^1(\Omega) \subset L^2(\Omega).
\end{equation}
If we consider the energy functionals $\Phi_{L^2}(u)$ (i.e., $\Phi_{\cH^{0}}(u)$) and  $\Phi_{\cH^{-1}}(u)$ defined in~\eqref{eq:Hs_obj_old}, we have 
\begin{align*}
&\nabla_{L^2} \Big( \Phi_{L^2}(u) \Big)= \mathcal{A}^*(\mathcal{A}u - f_\sigma),\quad  \nabla_{H^1} \Big( \Phi_{L^2}(u)  \Big)= (I-\Delta)^{-1}\mathcal{A}^*(\mathcal{A}u - f_\sigma), \\  &\nabla_{L^2} \Big( \Phi_{\cH^{-1}}(u)  \Big)= \mathcal{A}^*(I-\Delta)^{-1}(\mathcal{A}u - f_\sigma).
\end{align*}
Correspondingly, we have the following three gradient flow equations: 
\begin{eqnarray}
u_t &=& -\mathcal{A}^*(\mathcal{A}u - f_\sigma)\hspace{2.1cm} (\text{$L^2$ gradient flow of $\Phi_{L^2}(u)$}), \label{eq:L2_L2}\\
u_t &=& -(I-\Delta)^{-1} \mathcal{A}^*(\mathcal{A}u - f_\sigma) \quad (\text{$\cH^{1}$ gradient flow of $\Phi_{L^2}(u)$}), \label{eq:H1_L2}\\
u_t &=& -\mathcal{A}^* (I-\Delta)^{-1}(\mathcal{A}u - f_\sigma)\quad (\text{$L^2$ gradient flow of $\Phi_{\cH^{-1}}(u)$}). \label{eq:L2_Hm1}
\end{eqnarray}
If $\mathcal{A}^*$ shares the same set of eigenfunctions as the Laplace operator $\Delta$, then $\mathcal{A}^* (I-\Delta)^{-1} = (I-\Delta)^{-1} \mathcal{A}^*$, and hence~\eqref{eq:H1_L2} is exactly equivalent to~\eqref{eq:L2_Hm1}. Even if $\mathcal{A}^*$ does not commute with $(I-\Delta)^{-1}$, one can still view $(I-\Delta)^{-1}$ as a smoothing (integral) preconditioning operator upon the residual $\mathcal{A}u - f_\sigma$, which we wish to reduce to zero no matter the objective function is $\Phi_{L^2}(u)$ or $\Phi_{\cH^{-1}}(u)$. To sum up, \eqref{eq:H1_L2} and~\eqref{eq:L2_Hm1} are similar in nature in terms of the spectral bias of the resulting gradient descent dynamics, which demonstrates the equivalence between the change of the gradient flow and the change of the objective function under certain circumstances. In contrast to~\eqref{eq:L2_L2}, both \eqref{eq:H1_L2} and~\eqref{eq:L2_Hm1} are equipped with the smoothing property due to the additional $(I-\Delta)^{-1}$ operator. 

\subsection{Changing the Rate of Convergence}\label{sect:conv-rate}
So far, our analysis has been focusing on how the $H^s$ norm is related to the data noise $n_\sigma$ and its regularization effects during the optimization process. In this section, we address another interesting property of the $H^s$ norm as the objective function: it may improve the rate of convergence in gradient descent.

Extending the $L^2$ gradient flow~\eqref{eq:L2_Hm1} to a general $\Phi_{H^s}(u)$ energy functional, we obtain a gradient flow equation with respect to $u$:
\begin{equation} \label{eq:L2_Hs}
    u_t = -\mathcal{A}^* \mathcal{P}_s^*\mathcal{P}_s(\mathcal{A}u - f_\sigma),
\end{equation}
where $\mathcal{P}_s = (I-\Delta)^{s/2}$. Minimizing the $\Phi_{H^s}(u)$ energy functional~\eqref{eq:Hs_obj_old} is equivalent to reducing the $H^s$ norm of the residual $\mathcal{R}:=\mathcal{A}u - f_\sigma$. Based on~\eqref{eq:L2_Hs}, we have that
\begin{equation} \label{eq:L2_Hs_R}
    \mathcal{R}_t =  \mathcal{A} u_t =  - \mathcal{A} \mathcal{A}^* \mathcal{P}_s^*\mathcal{P}_s \mathcal{R}.
\end{equation}
The decay rate of the residual $\mathcal{R}$ is directly determined by the spectral property of the linear operator $\mathcal{A} \mathcal{A}^* \mathcal{P}_s^*\mathcal{P}_s$. After discretization, \eqref{eq:L2_Hs_R} becomes
\begin{equation*} 
  R^{(k)} =  (I - \eta E_s) R^{(k-1)} =  (I - \eta E_s)^k R^{(0)}, \quad E_s = A A^\top  P_s^\top P_s,
\end{equation*}
where $I$ is the identity matrix and $\eta$ is a properly chosen step size in gradient descent. As a result,
\begin{equation*} 
  \|R^{(k)}\|_2 =  \| (I - \eta E_s)^k R^{(0)} \|_2 \leq (1-\eta \lambda_{\text{min}})^k \|R^{(0)} \|_2,
\end{equation*}
where $\lambda_{\text{min}}$ is the minimum eigenvalue of $E_s$, which consequently depends on the choice of $s$. Given a fixed forward operator $\mathcal{A}$, by properly choosing $s$, we may improve the convergence rate by increasing $\lambda_{\text{min}}$. For example, if $\mathcal{A}u = \Delta u$, choosing the $H^{-2}$ norm as the objective function yields the fastest convergence among the class of $H^s$ norms~\cite{yang2020anderson}.

\section{Numerical Computation of the $H^s$ Norms}\label{sec:Hs_dist}
In this section, we present three numerical methods for computing the general $H^s$ norms of any $s\in \R$. The first one (in~\Cref{sect:firstHs}) applies to periodic functions defined on a domain, which is either the entire $\R^d$ or  a compact subset of $\R^d$, denoted by $\Omega$. We are mainly interested in periodic functions to align with a fast implementation of convolution that assumes the periodic boundary condition. 
In addition, we discuss the functions with zero Neumann boundary condition in~\Cref{sect:secondHs} and integer-valued $s$  in~\Cref{sect:thirdHs}. 


\subsection{Through the Discrete Fourier Transform}\label{sect:firstHs}
Recall that the Hilbert space $\cH^s(\mathbb{R}^n)$, $s\in\R$, is equipped with the norm~\eqref{eq:Hs_def}. If we compute the $H^s$ norm of a periodic function $f \in \cH^s$ defined on the entire $\R^d$, or equivalently, defined on $\Omega\subset \R^d$, we have
\begin{equation}\label{eq:Hs_periodic}
\left\| f\right\|_{\cH^{s}(\mathbb{R}^n)} = \left\| \mathcal{P}_s f\right\|_{L^2(\mathbb{R}^n)} \approx \left\|  {P}_s f\right\|_{L^2(\mathbb{R}^n)}, 
\end{equation}
where $\mathcal{P}_sf = \mathcal{F}^{-1}\left[(1+|\xi|^2)^{s/2}\mathcal{F}f\right]$ and ``$\approx$'' indicates the approximation by discretization. The discretization of the linear operator $\mathcal{P}_s$, denoted as $P_s$, can be computed explicitly through diagonalization, or implicitly, through the fast Fourier transform. For the former, the discretization of $\mathcal{F}$ is the discrete Fourier transform (DFT) matrix, while the discretization of $\mathcal{F}^{-1}$ is its conjugate transpose. The discretization of $(1+|\xi|^2)^{s/2}$ is correspondingly a diagonal matrix.


\subsection{Through the Discrete Cosine Transform}\label{sect:secondHs}
If we are interested in computing the $H^s$ norm of non-periodic functions on  the domain $\Omega$ that is a compact subset of $\R^d$, we adopt the zero Neumann boundary condition~\cite{schechter1960negative} as the boundary condition for the Laplacian operator.
As a result, rather than DFT, a consistent definition is through the discrete cosine transform (DCT) due to its relationship with the discrete Laplacian on a regular grid associated with the zero Neumann boundary condition, i.e., 
\begin{equation} \label{eq:Hs_DCT}
\|f\|_{\cH^{s}(\Omega)} \approx \| \widehat{P}_s f\|_{L^2(\Omega)},\quad \widehat{P}_s =  {C}^{-1} (I-\Lambda)^{s/2} {C},
\end{equation} 
where ${C}$ and ${C}^{-1}$ are matrices representing the DCT and its inverse, respectively, and $\Lambda$ is a diagonal matrix whose diagonal entries are eigenvalues of the  discrete Laplacian with the zero Neumann boundary condition. One may observe that~\eqref{eq:Hs_DCT} shares great similarity with~\eqref{eq:Hs_def} except for the facts that DFT is replaced with DCT and 
the diagonal matrix also varies according to eigenvectors and eigenvalues of the discrete Laplacian with different boundary conditions.

\subsection{Through Solving a Partial Differential Equation}\label{sect:thirdHs}
Let $\Omega \subset \mathbb{R}^n$ be a bounded Lipschitz-smooth domain. The Hilbert space $\cH^{s}(\Omega)$ is the same as the Sobolev space $W^{s,2}(\Omega)$ for all integers $s \in \mathbb{Z}$; see~\cite[Sec.~7]{arbogastmethods}, i.e.,
\[
W^{s,2}(\Omega) = \{f|_\Omega: f\in W^{s,2}(\R^d)\} =  \{f|_\Omega: f\in H^s(\R^d)\} = H^{s}(\Omega).
\]
Consequently, we can define an equivalent norm for functions in $H^{s}(\Omega)$ through $\|\cdot\|_{W^{s,2}(\Omega) }$, which  involves differential operators with the zero Neumann boundary conditions~\cite{schechter1960negative}. When $s\in \mathbb{N}$, the computation of the $W^{s,2}(\Omega)$ norm should follow its definition in~\cref{def:Wkp_def} while the differential operator involved should be handled with the zero Neumann boundary condition.
In this case, one explicit definition of $\|f\|_{\cH^{-s}(\Omega)}$ via the Laplace operator~\cite{schechter1960negative,yang2020anderson} is given by
\begin{equation} \label{eq:dual_Hs}
\|f\|_{\cH^{-s}(\Omega)} =  \|u\|_{\cH^{s}(\Omega)}, 
\end{equation} 
where $u(x)$ is the solution to the following partial differential equation with the zero Neumann boundary condition~\cite[Section 3]{schechter1960negative}, 
\begin{equation}\label{eq:dual_Hs_PDE}
    \begin{cases}
\mathfrak{L}^{s} u(x) =f(x), &x\in\Omega, \\
\nabla u\cdot {\bf n} = 0, &x\in \partial\Omega,
    \end{cases}
\end{equation}
for $
       \mathfrak{L}^{s} =  \sum\limits_{|\alpha|\leq s} (-1)^{|\alpha|} D^{2\alpha} .
$

We may define the operator $\mathfrak{L}^{-s}$  by setting $u = \mathfrak{L}^{-s} f$. Combining~\eqref{eq:dual_Hs} and~\eqref{eq:dual_Hs_PDE}, we have
\begin{equation}
    \|f\|^2_{\cH^{-s}(\Omega)} = \langle u, f \rangle_{L^2(\Omega)}   =   \langle \mathfrak{L}^{-s} f, f \rangle_{L^2(\Omega)}  =    \| \widetilde{\mathcal{P}}_s f\|^2_2,\quad \text{where}\quad  \widetilde{\mathcal{P}}_s^* \widetilde{\mathcal{P}}_s =  \mathfrak{L}^{-s}.
\end{equation}
We may also denote $\widetilde{\mathcal{P}}_s  = \mathfrak{L}^{-s/2}$. The numerical discretization of $\widetilde{\mathcal{P}}_s $ is denoted as $\widetilde{P}_s$.

Note that~\eqref{eq:Hs_def} and~\eqref{eq:dual_Hs} do not yield precisely the same norm given $f\in \cH^s(\R^d)$ with $s\in\mathbb{Z}$. For example, when $s=-2$ and $d=2$, the definition~\eqref{eq:Hs_def} depends on the integral operator $(I-\Delta)^{-1}$ based on the definition of the $\cH^{-s}(\Omega)$ norm, while the definition~\eqref{eq:dual_Hs} depends on the integral operator $(I-\Delta + \Delta^2)^{-1/2}$ based on the definition of the $W^{-s,2}(\Omega)$ norm in~\eqref{def:Wkp_def}. However, the leading terms in both definitions match. Thus, they are equivalent norms for functions that belong to the same  functional space $\cH^{s}(\Omega) = W^{s,2}(\Omega)$ given a fixed $s$. 
We remark that the $H^s$ norms with non-integer $s$ cannot be calculated through PDEs; instead, one should refer to Section~\ref{sect:secondHs}.




\section{Experiments}\label{sect:exp}

In this section, we first presents the denoising results of low-frequency noise arisen in geographical images in~\Cref{sect:exp_lf_noise}, followed by a nonlinear geophysical inverse problem in~\Cref{sect:exp_fwi}. In both examples, there is no regularization term in the objective function, so the implicit regularization effects purely come from the $H^s$ norm as the data fidelity term.

\begin{figure}
\centering
\subfloat[Noisy Input, PSNR= 27.53]{\includegraphics[width = 0.5\textwidth]{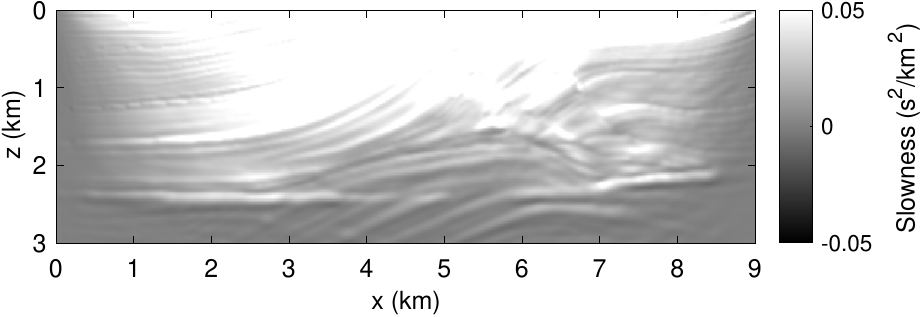}\label{fig:Marm_input}}
\subfloat[$\dot{\cH}^1$, PSNR= 36.34]{\includegraphics[width = 0.5\textwidth]{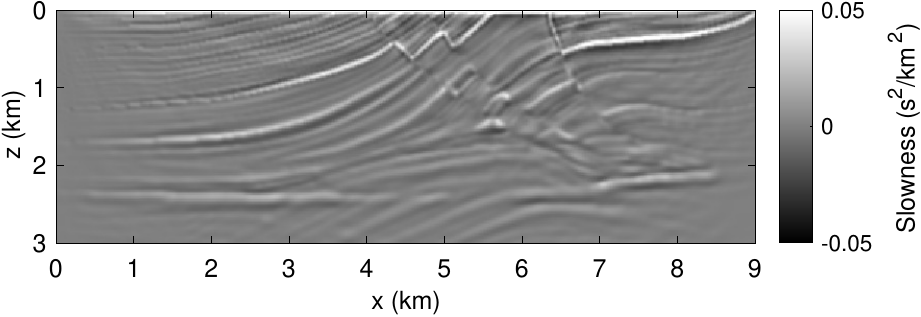}\label{fig:Marm_H1}}\\
\subfloat[$\dot{\cH}^2$, PSNR= 37.54]{\includegraphics[width = 0.5\textwidth]{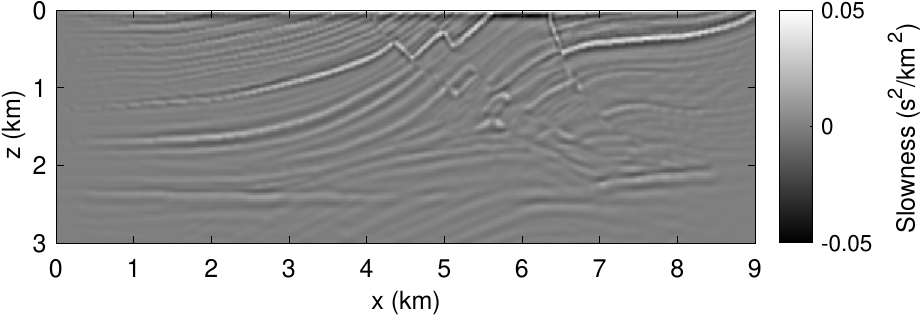}\label{fig:Marm_H2}}
\subfloat[$\dot{\cH}^3$, PSNR= 37.97]{\includegraphics[width = 0.5\textwidth]{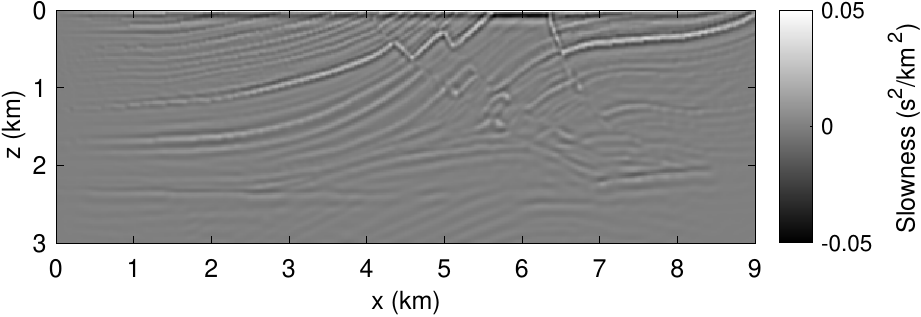}\label{fig:Marm_H3}}
\caption{Marmousi RTM image denoising using different $\dot{\cH}^s$ semi-norms as the data fidelity term.\label{fig:Marm_denoise}}
\end{figure}

\subsection{Geophysical Image  Denoising}\label{sect:exp_lf_noise}
We present a denoising example from a seismic application, in which the noise is mostly of low frequencies. Reverse-time migration (RTM)~\cite{claerbout1971toward} is a prestack two-way wave-equation migration to illustrate complex structure, especially strong contrast geological interfaces such as environments involving salts. Conventional RTM uses an imaging condition which is the zero time-lag cross-correlation between the source and the receiver wavefields. It overcomes the difficulties of ray theory and further improves image resolutions by replacing the semi-analytical solutions to the wave equation with fully numerical solutions for the full wavefield. 

However, artifacts are produced by the cross-correlation of source-receiver wavefields  propagating in the same direction. Specifically, migration artifacts appear at shallow depths, above strong reflectors, and severely mask the migrated structures; see~\Cref{fig:Marm_input}. They are generated by the cross-correlation of reflections, backscattered waves, head waves, and diving waves~\cite{zhang2009practical}. We are interested in reducing the strong low-frequency noise in the input data by minimizing the objective function~\eqref{eq:Hs_obj_semi}, where the linear operator $\mathcal{A}$ is the identity. Based on the discussion in Section~\ref{sec:Hs_analysis}, it is beneficial to use strong norms (i.e., $s>0$) to suppress the low-frequency noise. Here, we consider $\dot{\cH}^1$, $\dot{\cH}^2$ and $\dot{\cH}^3$ with the corresponding results shown in~\Crefrange{fig:Marm_H1}{fig:Marm_H3}, respectively. We  quantitatively measure the reconstruction performance  in terms of the peak signal-to-noise ratio (PSNR), which is defined by 
	\begin{equation*}
 \mbox{PSNR}( u^\ast, \tilde{u}) := 20 \log_{10} \frac{N M}{\| u^\ast-\tilde{ u}\|_2^2},
	\end{equation*} 
	where $u^\ast$ is the restored image, $\tilde{ u}$ is the ground truth, and $N, \  M$ are the number of pixels and the maximum peak value of $\tilde{ u},$ respectively. According to PSNR, using the $\dot{\cH}^3$ norm as the objective function produces the best recovery. We also demonstrate that all the three strong semi-norms can effectively
suppress the low-frequency noise in~\Cref{fig:Marm_input} without changing the  reflecting features of the underlying image.

\subsection{Full Waveform Inversion}\label{sect:exp_fwi}
Here we present a full waveform inversion (FWI) example. It is a nonlinear inverse problem where one aims to invert parameter $u$ (usually the wave velocity) given the observed data $g$ (usually the wave pressure field) through a nonlinear relationship $\mathcal{F}(u) = g$. The forward operator $\mathcal{F}$ is implicitly given through the wave equation constraint. For a more detailed introduction of this inverse problem, we refer to~\cite{virieux2009overview}. 

The nonlinear inverse problem is often reformulated as a PDE-constrained optimization problem where one aims to find the optimal $u$ by minimizing the difference between the observed data $g$ and the simulated data $\mathcal{F}(u)$ evaluated at the current prediction of $u$. While the least-squares method, i.e., using the squared $L^2$ norm as the data fidelity term, has been the conventional choice, and an additional regularization term is often added, here we consider only the general $H^s$ data fitting term as the objective function. That is, 
\begin{equation}\label{eq:FWI}
   \min_{u} \frac{1}{2}\| \mathcal{F}(u) - g\|^2_{H^s}.
\end{equation}
We perform optimization using different $s$ values and demonstrate its impacts on the inversion.
 
\begin{figure}
     \centering
     \subfloat[true velocity]{\includegraphics[width = 0.49\textwidth]{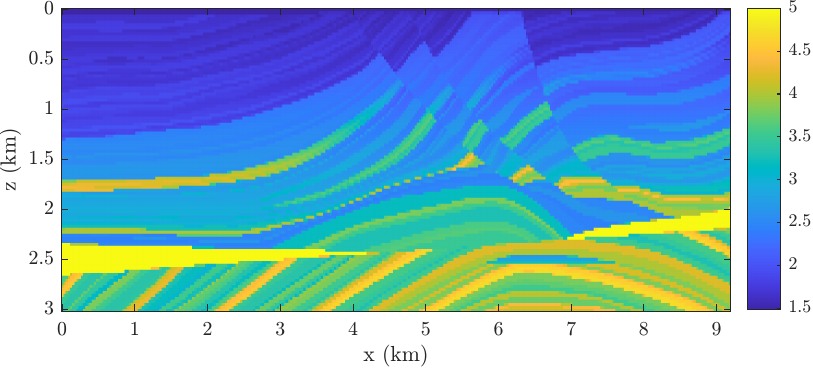}\label{fig:fwi_true}}
     \subfloat[initial guess]{\includegraphics[width = 0.49\textwidth]{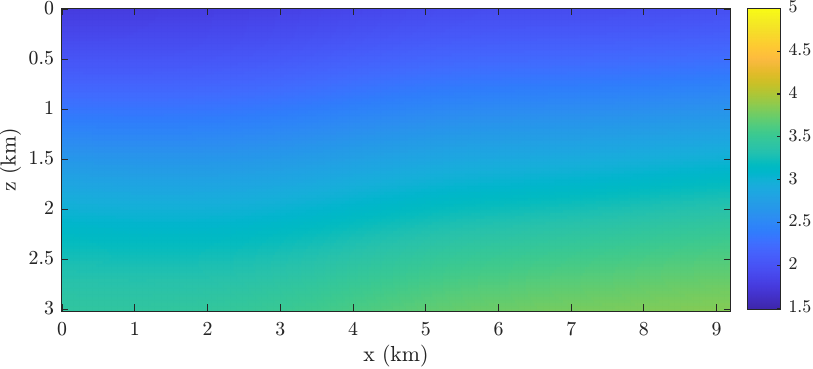}\label{fig:fwi_init}}\\
     \subfloat[$s=0$]{\includegraphics[width = 0.49\textwidth]{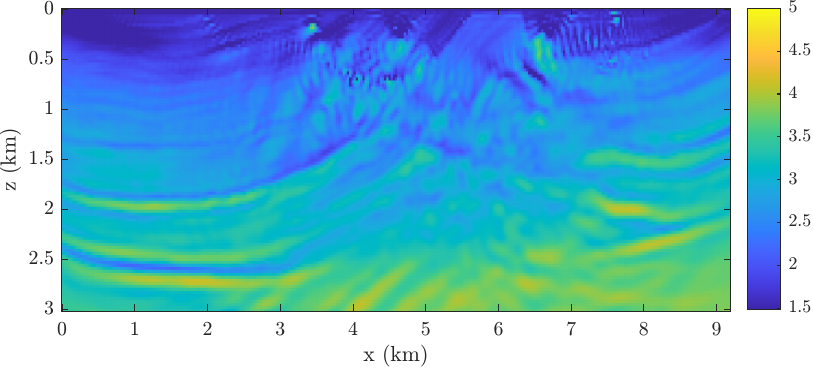}\label{fig:fwi_s0}}
     \subfloat[$s=-0.5$]{\includegraphics[width = 0.49\textwidth]{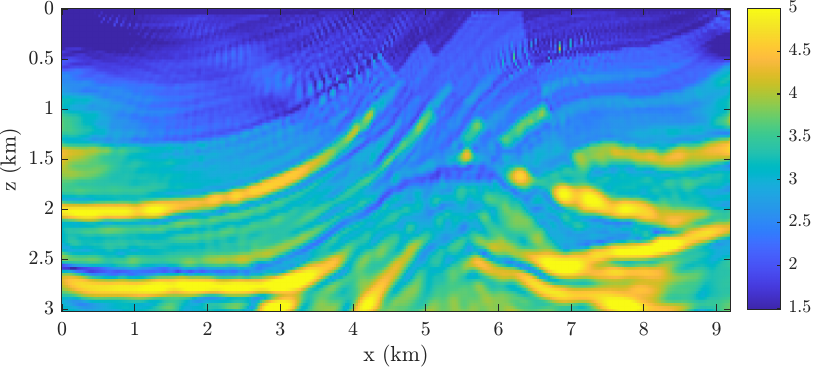}\label{fig:fwi_s05}}\\
     \subfloat[$s=-1$]{\includegraphics[width = 0.49\textwidth]{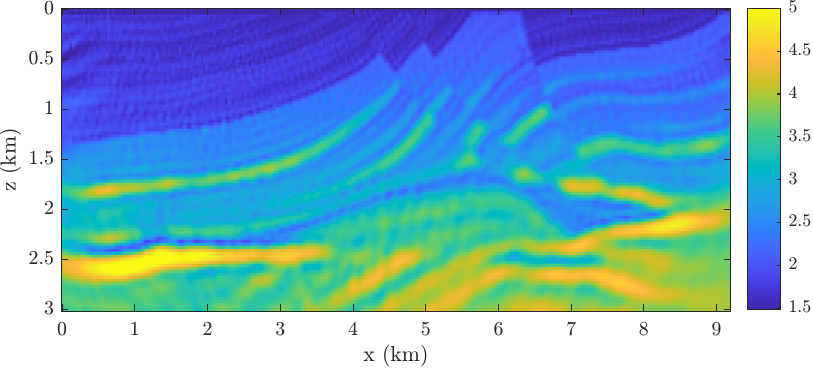}\label{fig:fwi_s1}}
     \subfloat[$s=-1$ then $s=0$]{\includegraphics[width = 0.49\textwidth]{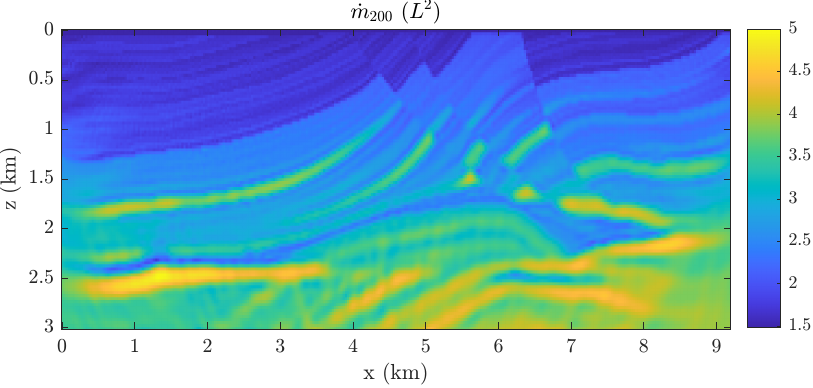}\label{fig:fwi_s1s0}}
     \caption{FWI example for Section~\ref{sect:exp_fwi}: (a)~true velocity; (b)~initial guess; (c)-(e)~reconstructions after 200 iterations using the $L^2$, $H^{-0.5}$, and the $H^{-1}$ norms, respectively; (f) reconstruction using the $H^{-1}$ in the first $100$ iterations followed by another $100$ iterations using the $L^2$ norm.\label{fig:fwi_all}}
\end{figure}

The true velocity parameter is presented in~\Cref{fig:fwi_true} and all the tests start with the same initial guess shown in~\Cref{fig:fwi_init}. We use the L-BFGS method~\cite{nocedal2006numerical} to solve for \eqref{eq:FWI} and manually stop the iterative process after $200$ iterations. The inversion result using the $L^2$ norm (corresponding to $s=0$) is shown in~\Cref{fig:fwi_s0}. It converges to a local minimum with many wrong features compared to the ground truth. Similarly, when using the $H^{-0.5}$ norm, the layers in the recovered subsurface image in~\Cref{fig:fwi_s05} do not match their true locations, despite a slight improvement from the $L^2$-based result. When using the $H^{-1}$ norm, the reconstruction is qualitatively much better as the structural properties of the inverted velocity image become very close to the ground truth, as one can see in~\Cref{fig:fwi_s1}. 

Since this is a nonlinear inverse problem, the resulting optimization problem~\eqref{eq:FWI} is highly nonconvex. The problem that the iterates are trapped at the local minima 
is often referred to as cycle skipping in FWI~\cite{virieux2009overview}. We expect that the change of the objective function modifies the optimization landscape.  It is well-known that low-frequency components of the wave data are less likely to suffer from cycle skipping~\cite{bunks1995multiscale}. As we have discussed in~Section~\ref{sec:Hs_analysis}, when $s<0$, the $H^s$ norm has a natural bias towards the low-frequency content of the input, and the smaller the $s$, the stronger the bias. Hence, it is not surprising to see that with the same initial guess, $H^{-1}$ norm as the objective function can converge to the global minima while the $L^2$ norm and the $H^{-0.5}$ norm get stuck at local minima. 

On the other hand, the $H^{-1}$ inversion in~\Cref{fig:fwi_s1} lacks high resolution despite having most of the correct features. Again, it is a property of the weak norm $(s<0)$. It is usually the high-frequency components of the data $g$ that resolve the sharp features in the reconstructed parameter $u$. However, the high-frequency components of the data, including both the useful physical information and the high-frequency noise, are given much smaller weight in a weak norm, resulting in a low-resolution reconstruction. We have commented on this phenomenon earlier in~\Cref{rmk:goodbad}. This dilemma can be mitigated by performing a transition of the objective function. For example, one can first use the $H^{-1}$ norm as the objective function to take advantage of the bigger basin of attraction. Once the iterate is close to the ground truth, one can switch to stronger norms such as the $L^2$. In~\Cref{fig:fwi_s1s0}, we perform a transition of the objective function from $H^{-1}$ after $100$ iterations of L-BFGS to the $L^2$ norm for another $100$ iterations. The resolution of the reconstruction is visibly improved compared to $200$ iterations of the $H^{-1}$ norm alone as shown in~\Cref{fig:fwi_s1}. A more rigorous analysis on how to adaptively update $s$ will be left to future work.

\section{A Case Study of Using Total Variation}
\label{sect:TVHs}

The natural implicit regularization effects of the $H^s$ norm could be further enhanced by combining with a regularization term, such as the TV regularization. There have been two main directions related to the combination in the literature.

First, minimizing the total variation energy under the general $H^s$ Sobolev space has been studied both numerically and theoretically~\cite{giga2010very,giga2017duality,kim2009image,schonlieb2015partial}. In such frameworks, the objective function is solely the TV energy, while the model parameter $u$ is assumed to belong to the $H^s$ functional space. Our work here is different from the literature since we fix the parameter space to be $L^2$ and consider the objective function to be $H^s$ or possibly $H^s$ together with a regularization term. As a result, the objective function explicitly includes the $H^s$ norm, equivalent to the assumption that the data space (as opposed to the model parameter space) is $H^s$. 

The second main direction in the literature is more relevant to our work. Combining the $H^{-1}$ data fitting term with the TV regularization was first studied in~\cite{osher2003image} and later generalized to any negative Sobolev norm in~\cite{lieu2008image}. The literature mainly focuses on image decomposition by using  TV to single out a cartoon (piece-wise constant) image and the $H^s$ norms for oscillatory components like textures and noises. Two recent works~\cite{huska2021variational,cicone2022jot} further propose to decompose a signal or an image into three components: a piece-wise constant component, a smooth (low-oscillating) component, and a high oscillatory component, the last of which is modeled by $H^{-1}$.

We advocate using the data fidelity term of the $H^s$ norm by itself as an implicit regularization effect on images. However, the frequency biases induced by the $H^s$ norm do not work so well on natural images due to complicated structures that spread out the entire frequency domain. As a result, image restoration requires an explicit regularization term to ensure satisfactory results. To this end, we present a proof-of-concept idea by incorporating the TV regularization together with the $H^s$ based data fidelity term. As $H^s$ reduces to the $L^2$ metric for $s=0$, we expect any regularization term combined with the $H^s$ would outperform the one with the standard least-squares misfit by treating $s$ as a tunable hyperparameter.

We also present a new algorithm to minimize the $H^s$ norm with the TV regularization based on ADMM, 
as detailed in~\Cref{sect:Hs+TV}. Under this efficient algorithmic framework, we then numerically investigate the power of combining the $H^s$ data-fitting term together with the TV regularization by presenting the deblurring examples in~\Cref{sect:deblurring}. The numerical results demonstrate that $H^s$+TV, as a more general framework, outperforms the traditional $L^2$+TV, making it a promising choice in image processing.



\subsection{An Numerical Algorithm for Minimizing TV regularization and $H^s$ Data Fitting Term}
\label{sect:Hs+TV} 


We revisit the celebrated TV regularization~\cite{rudin1992nonlinear} for image restoration that minimizes the following energy functional,
\begin{equation}\label{eq:image_model}
J(u) = 	 \frac{\lambda}{2} \|\mathcal{A} u -f_\sigma\|_{H^s}^2 + \mu \| \nabla u \|_1,
\end{equation}
where $\lambda,\mu\in \R^+$ are scalars balancing the data fitting term and the regularization term. We include two parameters $\lambda, \mu$ for the ease of disabling either one of them in experiments. We consider that the linear operator  $\mathcal{A}$ is either the identity operator for the denoising task or a convolution operator for the deblurring task, and $f_\sigma$ is the noisy (blurry) data. Osher, Sol\'e and Vese first proposed the framework~\eqref{eq:image_model} for the case $s=-1$~\cite{osher2003image}, which was later generalized by Lieu and Vess in~\cite{lieu2008image} to any $s<0$. Here, we extend the framework and apply it to any $s\in\mathbb{R}$. Moreover, we regard $s$ as a tunable hyperparameter, together with $\lambda$ and $\mu$ in~\eqref{eq:image_model}.

 We discuss the discretization of the model \eqref{eq:image_model}.  
	Suppose a two-dimensional (2D) image is defined on an $m\times n$ Cartesian grid. By using a standard linear index, we can represent a 2D image as a vector, i.e., the $((i-1)m+j)$-th component  denotes the intensity value at pixel $(i,j).$  	We define a discrete gradient operator,
	\begin{equation}\label{eq:gradient}
		\mathbf{D} u:= \left[\begin{array}{l}
	D_x\\
	D_y
		\end{array}
		\right]  u,
	\end{equation}
where $D_x,D_y$ are the finite forward difference operator with the periodic boundary condition in the  horizontal and vertical directions, respectively.
We adopt the periodic boundary condition for finite difference scheme to align with the periodic boundary condition when implementing the discrete convolution operator $A$ by the fast Fourier transform (FFT). 
We denote  $N := mn$ and the Euclidean spaces by  $\mathcal X:=\mathbb{R}^{N}, \mathcal Y:=\mathbb{R}^{2N}$, then $u\in \mathcal X,$ $Au\in \mathcal X,$ and $\mathbf{D} u\in \mathcal Y$.


The $\cH^s$ norm can be expressed in terms of the weighted norm, which is equivalent to the multiplication of $\mathbf{P}_s$, the discrete representation of the operator $\mathcal{P}_s$. Given the choice of $s$ and the particular boundary condition, 
we can select a preferable way of implementing 
$\mathbf{P}_s$ as any of the three types of matrices ${P}_s$, $\widehat {P}_s$, and $\widetilde {P}_s$  discussed in Section~\ref{sec:Hs_dist}. 
To align with the periodic boundary condition used for $\mathbf D$ and $A$, we choose 
 $\mathbf{P}_s = P_s$.
In summary, we  obtain the following objective function in a discrete form,
\begin{equation}\label{eq:Hs_obj}
J(u) = 	 \frac{\lambda}{2} \| \mathbf{P}_s(A u -f_\sigma)\|_{2}^2 + \mu \| \mathbf{D} u \|_1.
\end{equation}

There are a number of optimization algorithms available to minimize $J(u)$ in order to find the optimal solution $u$, such as the Newton's method, the conjugate gradient descent method, and various quasi-Newton methods~\cite{esser2010general,goldstein2009split,nocedal2006numerical}. Here,  we present the alternating direction method of multipliers (ADMM)~\cite{boyd2011distributed,glowinski1975approximation},
by introducing an auxiliary variable $d$ and studying an equivalent form of~\eqref{eq:Hs_obj}
\begin{equation}\label{equ:split_model_uncon}
	\min_{u \in\mathcal X,  d\in\mathcal Y} \quad \mu \|  d \|_1+\frac{\lambda}{2} \|\mathbf{P}_s(A u -f_\sigma)\|_2^2 \quad \mathrm{s.t.} \quad  d = \mathbf{D} u.
\end{equation}
The corresponding augmented Lagrangian function is expressed as
\begin{equation}\label{eq:AL4L1uncon}
	\mathcal{L}(u, d; v) = \mu\| d \|_1+\frac{\lambda}{2} \|\mathbf{P}_s(A u -f_\sigma)\|_2^2+\langle \rho  v,\mathbf{D} u - d\rangle + \frac{\rho}{2}\| d - \mathbf{D} u \|_2^2,
\end{equation}
with a dual variable $ v$ and a positive parameter $\rho.$ The ADMM framework involves the following iterations,
\begin{equation} \label{ADMML1_uncon}
	\left\{\begin{array}{l}
	u^{(k+1)}=\arg\min_u \mathcal{L}(u,  d^{(k)};  v^{(k)}),\\
	 d^{(k+1)}=\arg\min_{ d} \mathcal{L}(u^{(k+1)},  d;  v^{(k)}),\\
	 v^{(k+1)} =   v^{(k)} + \mathbf{D}  u^{(k+1)} -  d^{(k+1)}.
	\end{array}\right.
\end{equation}

By taking the derivative of $\mathcal{L}$ with respect to $u$, we obtain a closed-form solution of the $u$-subproblem in \eqref{ADMML1_uncon}, i.e.,
\begin{equation}\label{ADMM_l1con_u}
	u^{(k+1)} = \left(\lambda A^T \mathbf{P}_s^T \mathbf{P}_s A + \rho \mathbf{D}^T \mathbf{D} \right)^{-1}\left(\lambda A^T \mathbf{P}_s^T \mathbf{P}_s f_\sigma  +  \mathbf{D}^T \big( d^{(k)} -\rho  v^{(k)} \big)\right).
\end{equation}
We remark that $-\mathbf{D}^T \mathbf{D}$ is the discrete Laplacian operator with the periodic boundary condition. In this case, the discrete operators (matrices), $A$, $A^T$, $\mathbf{P}_s^T \mathbf{P}_s$ and $\mathbf{D}^T \mathbf{D}$ all have the discrete Fourier modes as eigenvectors. As a result, the matrix $\lambda A^T \mathbf{P}_s^T \mathbf{P}_s A + \rho \mathbf{D}^T \mathbf{D}$ in~\eqref{ADMM_l1con_u} shares the Fourier modes as eigenvectors, and its inverse can be computed efficiently by FFT.


The $d$-subproblem in \eqref{ADMML1_uncon} has also a closed-form solution given by
\begin{equation}\label{ADMM_l1con_d}
	\h d^{(k+1)} = \mathbf{ shrink}\left(\nabla u^{(k+1)} + \h v^{(k)}, \frac{\mu}{\rho}\right),
\end{equation}
where 
$
\mathbf{shrink}(\h v, \beta) = \mathrm{sign}(\h v)\circ \max\left\{|\h v|-\beta, 0\right\}
$ with the Hadamard (elementwise) product
 $\circ$. 
 Finally, $ v^{(k+1)}$ is updated based on $u^{(k+1)}$  and $ d^{(k+1)}$. The iterative process continues until reaching the stopping criteria or the maximum number of iterations.

\subsection{Image Deblurring}\label{sect:deblurring}
We start this subsection by first expanding the deblurring example in~\Cref{sec:Hs_analysis}. In particular, we conduct a comprehensive study of the $H^s$ norms  with different choices of $s$ under a variety of noise levels and whether the TV regularization term is included in the objective function or not. We remark that the noise here is high-frequency Gaussian noise. The PSNR values in different settings of deblurring the Square image are recorded in~\Cref{tab:Hs_Compare}.

The first row of~\Cref{tab:Hs_Compare} is about  
 the reconstruction without using TV from noise-free data, i.e., $\sigma = 0$. 
All the PSNR values  
 are all over 190, which implies the perfect recovery (subject to  numerical round-off errors). 
In this noise-free case, the reconstruction is a standard (weighted) least-squares solution. Furthermore, the choice of the data-fitting term does not   affect the minimizer of the optimization problem, though the convergence rate may differ. As seen in~\Cref{fig:square_deblur_nonoise_different_s}, the same number of gradient descent iterations yields different sharpness  when $s$ varies. 

\begin{figure}
    \centering
    \includegraphics[width = 1.0\textwidth]{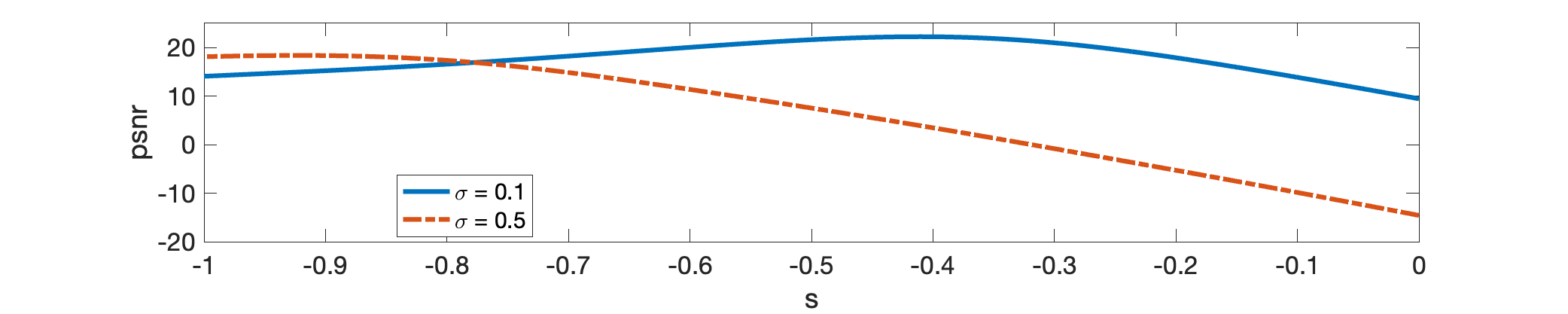}
    \caption{Illustrating how the PSNR value depends on different $\cH^s$ norms for deblurring the Square image without regularization. The optimal $s$ varies with the noise intensity. For a larger noise variance, it is preferable to select a weaker norm (corresponding to a smaller $s$).}
    \label{fig:square_deblur_diff_s}
\end{figure}

Still without the regularization term, we examine the denoising results using the noisy blurry data and record the PSNR values in the second and the third rows of~\Cref{tab:Hs_Compare}. These quantitative values reflect that the reconstruction results after a fixed number of gradient descent iterations \eqref{eq:Hs_GD} differ drastically with respect to different $s$ values, as also illustrated in~\Cref{fig:square_deblur}.
We plot the PSNR values with more $s$ values in~\Cref{fig:square_deblur_diff_s} than those documented in~\Cref{tab:Hs_Compare}, which further illustrates that the optimal choice of $s$ depends on the noise level. 


The effect of the TV regularization is presented in the last two rows of~\Cref{tab:Hs_Compare}. On one hand, TV significantly improves the results over the model without TV. 
On the other hand, using the optimal $\cH^s$ norm as the data-fitting term together with TV outperforms the classic TV with the $L^2$ norm, as the former has an extra degree of freedom.

\setlength{\tabcolsep}{6pt}
\begin{table*}[t]
\centering
\begin{tabular}{|c|c|c|c|c|c|c|c|}
\hline
Add $TV$ 	& Noise $\sigma$	  & input & $s = 0$ &  $s =-0.25 $ & $s =-0.5 $ & $s =-0.75 $ & $s =-1$\\
\hline \hline
No &   0 &    24.25 &  194.57 & 194.57 & 194.57 & 194.57 & 194.57\\
No &   0.1 & 18.61& 9.46  & 19.62  & 21.63  & 17.38  & 14.09\\
No &   0.5 & 5.95 & -14.57 &  -3.05  &  7.54 &  16.29  & 18.12\\
Yes &   0.1 & 18.61&    39.03  & 39.49  & 39.85  & 40.16  & 40.39\\
Yes &   0.5 & 5.95 & 27.67  & 27.99 &  28.23  & 28.39 &  28.44\\
\hline
\end{tabular}
\caption{Deblurring the Square image comparison among different $\cH^s$ norms in terms of PSNR. Visual results corresponding to the second and the third rows are shown in \cref{fig:square_deblur}.} \label{tab:Hs_Compare}
\end{table*}

\begin{figure}[t]
\centering
\subfloat[Noisy Input]{\includegraphics[width = 0.25\textwidth]{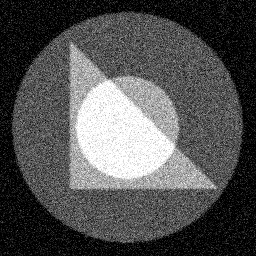}}
\hspace{0.1cm}
\subfloat[TV]{\includegraphics[width = 0.25\textwidth]{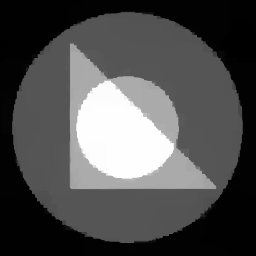}}
\hspace{0.1cm}
\subfloat[Hyper]{\includegraphics[width = 0.25\textwidth]{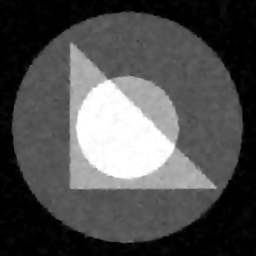}}\\
\subfloat[BM3D]{\includegraphics[width = 0.25\textwidth]{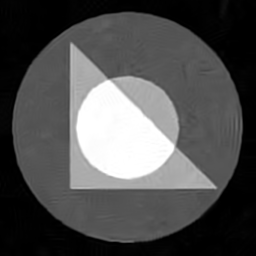}}
\hspace{0.1cm}
\subfloat[WAI]{\includegraphics[width = 0.25\textwidth]{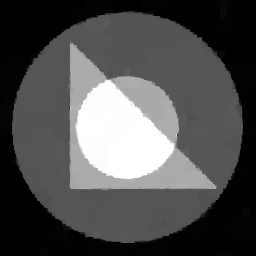}}
\hspace{0.1cm}
\subfloat[proposed]{\includegraphics[width = 0.25\textwidth]{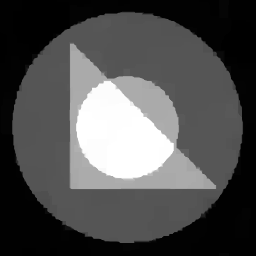}}\\
\caption{Comparison of deblurring the Circles image  with a $7\times 7$ Gaussian blur  and additive Gaussian noise of $\sigma = 0.1.$\label{fig:circle_deblur}}
\end{figure}

\begin{figure}[t]
\centering
\subfloat[Noisy Input]{\includegraphics[width = 0.25\textwidth]{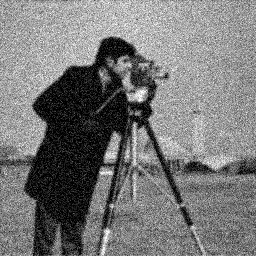}}
\hspace{0.1cm}
\subfloat[TV]{\includegraphics[width = 0.25\textwidth]{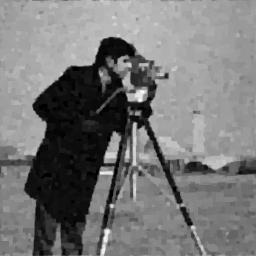}}
\hspace{0.1cm}
\subfloat[Hyper]{\includegraphics[width = 0.25\textwidth]{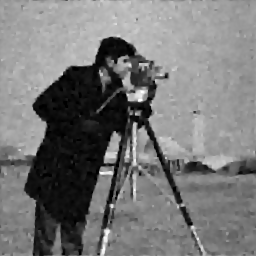}}\\
\subfloat[BM3D]{\includegraphics[width = 0.25\textwidth]{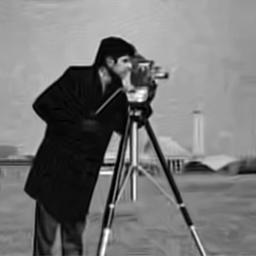}}
\hspace{0.1cm}
\subfloat[WAI]{\includegraphics[width = 0.25\textwidth]{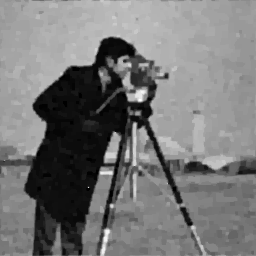}}
\hspace{0.1cm}
\subfloat[proposed]{\includegraphics[width = 0.25\textwidth]{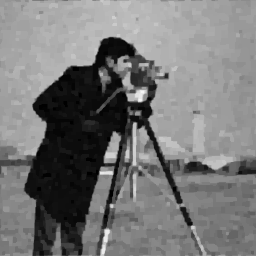}}
\caption{Comparison of deblurring the Cameraman image with $7\times 7$ Gaussian blur and additive Gaussian noise of $\sigma = 0.1.$\label{fig:cameraman_deblur}}
\end{figure}

We further test on two images: Circles and Cameraman, for image deblurring. The blurring kernel is fixed as a $7\times 7$ Gaussian function with the standard deviation of 1. 
By assuming the periodic boundary condition and using the Convolution Theorem, the linear operator $A$  can be implemented by FFT. We also consider two noise levels: $\sigma = 0.1$ and $0.2$ as the standard deviation of the additive Gaussian random noise. 
We compare the proposed approach $H^s$+TV with  TV, a hyper-Laplacian model (Hyper) \cite{krishnan2009fast}, a modification of BM3D from denoising to deblurring \cite{dabov2008image}, and a weighted anisotropic and isotropic (WAI) regularization proposed in~\cite{lou2015weighted}. We use the online codes of the competing methods: Hyper, BM3D, and WAI. For all the methods, we tune the parameters so that they can achieve the highest PSNR for each combination of testing image and noise level.  We record the PSNR values in \cref{tab:psnr_deblur} and present the visual results under a lower noise level ($\sigma=0.1$) in \cref{fig:circle_deblur,fig:cameraman_deblur}. The proposed approach works particularly well for images with simple geometries such as Circles, and is comparable to the state-of-the-art deblurring methods for the Cameraman image.

\setlength{\tabcolsep}{10pt}
\begin{table*}[t]
\centering
\begin{tabular}{ |c|c|c|c|c|c|c|c|}
\hline
Test image 	& $\sigma$	  & input &  TV &  Hyper & BM3D & WAI & proposed  \\
\hline \hline
\multirow{3}[6]*{\bf Circles}  & 0.1 &   19.78&   32.56&   30.61&   32.52&   31.96&   32.93
 \bigstrut\\
\cline{2-8}
& 0.2 &  13.91 &   29.84 &   28.10 &   29.97 &   29.78 &   30.03 \bigstrut\\
\hline \hline
\multirow{3}[5]*{\bf Cameraman}  & 0.1 &   18.96&   24.52&   24.54 &   25.49&   24.40&   24.53\bigstrut\\
\cline{2-8}
& 0.2 &  13.65&   22.89 &   22.75&   23.53 &   22.92 &   22.96
\bigstrut\\
\hline
\hline
\end{tabular}
\caption{Image deblurring comparison in terms of PSNR.} \label{tab:psnr_deblur}
\end{table*}

\section{Conclusions}\label{sect:conclusion}
In this paper, we proposed a novel idea of using the Sobolev ($H^s$) norms as a data fidelity term for imaging applications. We revealed implicit regularization effects offered by the proposed data fitting term rather than the commonly used regularization term. Specifically, we shall choose a weak norm ($s<0$) for high-frequency noises and a strong norm ($s>0$) for low-frequency noises. We discussed the connections between the Sobolev norm and the Sobolev gradient flow. From a Bayesian inference perspective, we analyzed the underlying noise assumption for a Sobolev norm as the data fidelity term. We further revealed that one could choose a proper Sobolev norm as an objective function to improve the convergence rate in gradient descent, achieving preconditioning effects. We presented three numerical schemes to compute the  $H^s$ norms under different domains and boundary conditions. 
Experimental results showed that the $H^s$ data fitting term alone as the objective function has implicit regularization effects on the performance of various inverse problems. Furthermore, the $H^s$ data fitting term combined with the TV regularization, i.e., $H^s$+TV, works particularly well for images with simple geometries and always outperforms the standard $L^2$+TV. In the framework of ADMM, one can efficiently minimize the $H^s$+TV model with  a tunable parameter $s$.

\bibliographystyle{siamplain}
\bibliography{references}

\end{document}